\documentclass[12pt]{amsart}

\usepackage{epsfig}
\usepackage{verbatim}
\usepackage{rotate}
\theoremstyle{definition}

\theoremstyle{remark}

\numberwithin{equation}{section}

\newcommand{\eps}{\varepsilon}

\newcommand{\der}[2]{\frac{\partial #1}{\partial #2}}
\newcommand{\derf}[2]{\frac{d #1}{d#2}}

\begin{document}

\title[]{Numerical investigation of a dipole type solution for unsteady groundwater flow with capillary retention and forced drainage}
\author{Eugene A. Ingerman \lowercase{\and} Helen Shvets \\
\emph{{D\lowercase{epartment of} M\lowercase{athematics},\\ U\lowercase{niversity of} C\lowercase{alifornia,} B\lowercase{erkeley and} \\
L\lowercase{awrence} B\lowercase{erkeley} L\lowercase{aboratory,} CA
94720 \\
\lowercase{eugening@math.berkeley.edu, shvets@math.berkeley.edu}}}}

\date{\today}
\dedicatory{}

\begin{abstract}
A model of unsteady filtration (seepage) in a porous medium with
capillary retention is considered.  It leads to a free boundary
problem for a generalized porous medium equation where the location of
the boundary of the water mound is determined as part of the
solution. The numerical solution of the free boundary problem is shown to
possess self-similar intermediate asymptotics.
On the other hand, the asymptotic solution can be obtained from a
non-linear boundary value problem.  Numerical solution of the
resulting eigenvalue problem agrees with the solution of the partial
differential equation for intermediate times. In the second part of
the work, we consider the problem of control of the water mound
extension by a forced drainage.
\end{abstract}

\maketitle
\pagestyle{plain}
\section{Introduction.}
In the present work two problems from the theory of filtration
through a horizontal porous stratum are considered. First we study a
short, but intense, flooding followed by natural outflow through the
vertical face of an aquifer. Further, we consider the possibility to
control the spreading of the water mound by use of forced drainage at
the boundary.

An important practical example of such a problem is groundwater mound
formation and extension following a flood, after a breakthrough of a
dam, when water (possibly contaminated) enters and then slowly extends
into a river bank.

\begin{figure}[htbp]
\centering \epsfig{figure=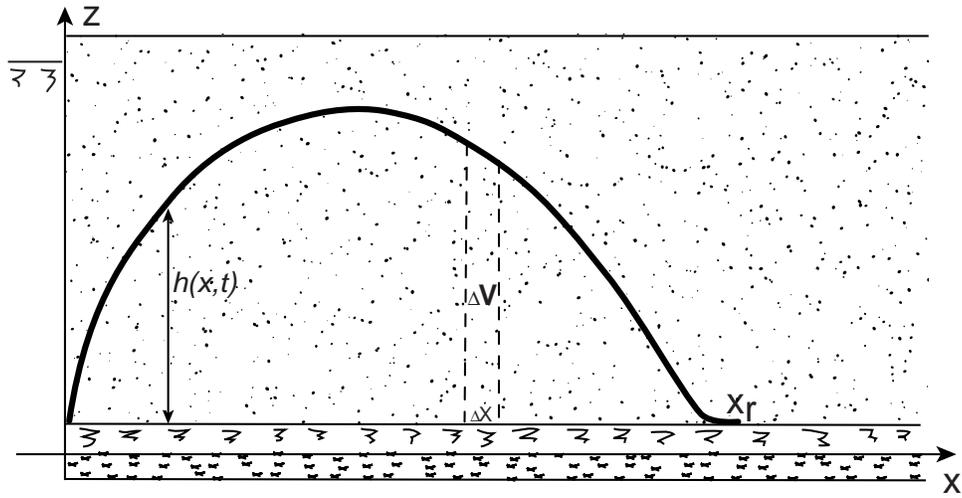,width=5in}
\caption{Groundwater dome extension in a porous medium.} \label{fig:fig1}
\end{figure}

Consider an aquifer that consists of a long porous stratum with an
impermeable bed at the bottom and a permeable vertical face on one
side {(Fig. \ref{fig:fig1})}. The space coordinate $x$ is directed
along the horizontal axis with $x=0$ at the vertical face.  A water
reservoir is located in the region $x<0$. We assume that the flow is
homogeneous in the y-direction. The height of the resulting mound is
denoted by $z=h(x,t)$.  The initial level of water in the stratum is
assumed to be negligible.

The problem is formulated as follows. At some time $t=-\tau<0$, the water
level at the wall begins to rise rapidly, and water enters the porous
medium.  By time $t=0$, the water level at the vertical face returns
to the initial one.We assume that the distribution at time $t=0$ is
given by $h(x,0)=h_0(x)$ and is concentrated over a finite region
$[0,d]$ (compactly supported). We also assume that $h_0(x)$ is concave
down and $h_0^2(x)$ is gently sloping.

In problem 1, water naturally seeps through the boundary back into the
reservoir, giving the boundary condition $h(0,t)=0$. Inclusion of the
effects of capillary retention into the model distinguishes our case
from the well know dipole-type problem. The numerical and asymptotic
solutions for the source-type boundary conditions were obtained in~[6].  Most recently, the
dipole-type problem with capillary retention was studied numericaly and analytically, using
Lie-group techniques, by B. Wagner in~[7].

Analysis and numeric computations show that in the case of natural
outflow, the water mound
is not extinguished in finite time. The outflow rate cannot be further
increased by lowering the level at the boundary.  In problem 2, in
order to control the spreading of the water mound, forced drainage is
introduced. The problem formulation was proposed in~[3],
where a complete mathematical derivation and rigorous analysis can be
found. The forced drainage can be implemented, for instance, by drilling a
number of holes or horizontal wells near the impermeable bottom.
In this way, an additional discharge rate is created, and the fluid level
becomes zero on some interval $[0, x_l(t)]$.

These are certainly highly idealized problems, but their solutions
allow one to extract the qualitative properties and to check the numerical
methods in solving more realistic problems.

\section{Porous medium equation with capillary retention.}
In the case of seepage and gently sloping profile $h^2 (x,t)$ and in
the absence of capillary retention, the model of flow in a porous
stratum is described by the Boussinesq equation ([4] see
also~[2],[1]):
\begin{equation}\label{eq0}
\partial_t h=\kappa \partial_{xx} h^2.
\end{equation}
Here $\kappa={k \rho g}/{2 \mu m }$, $k$ is the permeability of the
medium, $m$ its porosity (the fraction of the volume in the stratum
which is occupied by the pores), $\rho$ the fluid density, $\mu$ its
dynamic viscosity, and $g$ the acceleration of gravity.  According to
the hydrostatic law, water pressure $p=\rho g(h-z)$.  Then, the
total head
$H=p+\rho g z = \rho g h$ is constant throughout the height of the
mound. Under the assumption of seepage and gently sloping profiles
$h^2$, Darcy law is used to obtain the relation for the total flux
$q=-\frac{k}{\mu} \nabla H \cdot h$.

Mathematical properties of the Boussinesq equation are well
known~[5]. An essential feature of this equation is the
finite speed of disturbance propagation given a finite (compactly
supported) initial distribution.  Another important feature of this
equation is the existence of special self-similar solutions. The
graphs of such a solution for any two times $t_0$ and $t_1$ are
related via a similarity transformation~[1].  The special
solutions, themselves corresponding to certain, sometimes artificial,
initial and boundary conditions, are important because they provide
intermediate asymptotics for a wide class of initial value
problems. For these problems, the details of the initial distribution
affect the solution only in the beginning; after some time, the
solution approaches a self-similar asymptotics. The Boussinesq
equation has been studied extensively and a number of self-similar
solutions, for different boundary conditions, have been constructed
([2],~[3]).

Following~[1],~[6], the Boussinesq equation can be
modified to incorporate the effects of capillary retention into the
model. If we exclude the possibility of water reentering the region
that was filled with water at some earlier time and assume that
initially the stratum is empty, we have the following situation: when
water enters a pore, it occupies the entire volume, allowed by active
porosity; when water leaves the pore, a fraction $\delta$ of the pore
volume remains occupied by the trapped water. We assume that $\delta$
is constant. Let us denote the initial active porosity by $m$. Then, when
water is entering previously unfilled pores, the effective porosity is
$m$; when water is leaving previously water-filled pores, the
effective porosity becomes $m (1-\delta)$. Hence, in the presence of
capillary retention, porosity depends on the sign of $\partial_t
h(x,t)$. Notice, that permeability can be assumed unaffected, as the
effect of capillary forces on permeability is significant only for
small and/or dead-end pores, whose contribution to the total flux, in
the first approximation, can be neglected.

The rate of change in the amount of water $\Delta V$ inside a volume
element (Fig. \ref{fig:fig1}) is equal to:
\begin{equation}
        \Delta V= \begin{cases} m \der{h}{t}\Delta x & \text{if
        $\der{h}{t}>0$ }, \\ m (1-\delta) \der{h}{t} \Delta x &\text{if
        $\der{h}{t}<0$}. \end{cases}
\end{equation}
On the other hand, the rate of change in the volume of water due to the
flux through
the faces of a volume element
${(x,{x+\Delta x})}$ is equal to
\begin{equation}  \Delta x \> \partial_x(\frac{k}{\mu} \partial_x H \cdot h)=\Delta x
 \> \frac{k \rho g}{2 \mu}\partial_{xx} (h^2).
\end{equation}
We denote $\kappa_1=\frac{k \rho g}{2 \mu m}$ and $\kappa_2=\frac{k
\rho g}{2 \mu m (1-\delta)}$. Then, using the continuity of flux (no
sources inside the water mound) and the balance of mass we obtain:
\begin{equation}\label{eq1}
        \der{h}{t}= \begin{cases} \kappa_1 \partial_{xx}(h^2) &
        \text{if $\der{h}{t} > 0$ },\\ \kappa_2 \partial_{xx}(h^2) &
        \text{ if $\der{h}{t}<0$}. \end{cases}
\end{equation}
This is a nonlinear parabolic partial differential equation with
discontinuous coefficients, also known as the generalized porous medium
equation~[2],~[6].

Continuity of the flux $q=-\frac{\rho g k}{\mu} h \cdot
\partial_{x}h$ implies that at the mound tip $x_r$, where mound
height is zero, the flux is also zero. For problem 1, these
considerations lead to the following initial and boundary conditions
to supplement equation (\ref{eq1}):
\begin{eqnarray} \label{eq2}
h(x,0) &=& h_0(x) \ge 0 \text{~~~(where } h_0(x)=0 \text{ for } x \ge
d \text{)}, \nonumber \\
h(x,t) &=& 0 \text{, }\partial_{x}h^2(x,t) = 0 \text{ at $x=x_r$}, \\
h(0,t) &=& 0. \nonumber
\end{eqnarray}
The second line in~(\ref{eq2}) corresponds to the free boundary
conditions on the right boundary, $x_r(t)$, which is unknown a priori.

It should be noted that for the solution of equation~(\ref{eq0})
(but not for~(\ref{eq1})) with boundary conditions~(\ref{eq2}) the
dipole moment is constant:
\begin{equation}
Q =\int\limits_{0}^{\infty}{x h(x,t) \,d x} = C.
\end{equation}

We call equation (\ref{eq1}) with boundary conditions (\ref{eq2}) a
dipole-type problem. A similar problem, for source type initial
and boundary conditions was considered in~[6], see
also~[1].

For problem 2, the boundary conditions are changed to include the
forced drainage condition. The discharge rate $q_0(t)$, which is a
quantity that should be specified, determines the boundary condition
at the left free boundary $x_l$.
\begin{eqnarray} \label{eq2_1}
h(x,0) &=& h_0(x) \ge 0 \text{~~~(where } h_0(x)=0 \text{ for } x \ge
d \text{)},\nonumber \\ 
h(x,t)&=& 0,~\partial_{x}h^2(x,t)= 0 \text{ at $x=x_r$},\\
h(x,t)&=& 0,~\partial_{x}h^2(x,t)= \frac{q_0(t)}{m \kappa}
\text{ at $x=x_l$ }. \nonumber
\end{eqnarray}
The second and third lines in (\ref{eq2_1}) define, respectively, the free
boundary condition on the right boundary and the forced drainage condition
on the left boundary. Equation (\ref{eq1}) together with boundary
conditions (\ref{eq2_1}) define problem~2.

\section{Dimensional analysis of problem 1.}

The parameters in the problem are $h,x,t,\kappa_1,\kappa_2$, $d$ - the
initial width of the water mound, and $Q=Q(0)$ - the initial dipole
moment. We can take the dimensions as follows: $[h]=H$, $[x]=L$,
$[t]=T$. Then from equation (\ref{eq1}) we have
$[\kappa]=\frac{L^2}{TH}$. For the remaining parameters $[d]=L$,
$[Q]=H \cdot L^2$. The dimensions for $h$ and $x$ are set to be
independent.  This can be done because the differential equation
(\ref{eq1}) is invariant with respect to the following group of
transformations:
\begin{equation}
        x'=\alpha x, \text{ } t'=\frac{\alpha^2}{\gamma} t, \text{ }
        h'=\gamma h.
\end{equation}
The invariance insures that we can scale the units of measurement for
$h$, while keeping the units for $x$ unchanged.

The following dimensionless quantities can be obtained from these
parameters: $$ \Pi_1=\frac{x}{(Q \kappa_1 t)^{1/4}},
\Pi_2=\frac{d}{(Q \kappa_1 t)^{1/4}},
\Pi_3=\frac{\kappa_1}{\kappa_2}, \text{and } \Pi=h(\frac{\kappa_1 t} {Q})^{1/2}.
$$ It follows that $\Pi=F(\Pi_1,\Pi_2,\Pi_3)$.

Since for large times, $t \gg \frac{\Delta x^4}{Q \kappa_1}$, the
parameter $\Pi_2 \ll 1$, it would seem natural to set
$\Pi=f(\Pi_1,\Pi_3)$, as in the case of $\kappa_1=\kappa_2$, and look for a
solution of the form:
\begin{equation}
h = (\frac{Q^2}{\kappa_1 t})^{1/2}f(z,\frac{k_1}{k_2}), 
\text{ where }z = \frac{x}{(Q \kappa_1 t)^{1/4}}.
\end{equation}

However, this leads to a contradiction when we consider an ordinary
differential equation obtained from (\ref{eq1}):
\begin{equation}\label{eq6}
 (2f+\derf{f}{z}z)=
\begin{cases}
        -4\derf{(f^2)}{z^2} & \text{if $2 f+\derf{f}{z}z < 0$},\\
        -4\frac{\kappa_1}{\kappa_2}\derf{(f^2)}{z^2} & \text{if $2
        f+\derf{f}{z}z>0$}.
\end{cases}
\end{equation}
Multiplying both sides by $z$ we obtain an equation in total
differentials, which is readily solved:
\begin{equation}\label{eq7}
 fz^2= \begin{cases} -4(z\frac{d f^2}{dz} - f^2) + C_1& \text{if
 $2f+\frac{df}{dz}z < 0$}, \\
 -4\frac{\kappa_1}{\kappa_2}(z\frac{df^2}{dz} - f^2) + C_2& \text{if
 $2f+\frac{df}{dz}z > 0$}.
\end{cases}
\end{equation}
Observe that near $z=z_r$, where the height of the mound vanishes, the
first equation holds. At $z=z_r$, $h$ vanishes along with the flux,
which is proportional to $\der{h^2}{x}$. From the first equation at
$z_r$, we obtain that $C_1=0$. Similarly, evaluating the second
expression at $z=0$, where $h=0$, we find that $C_2=0$. Next,
evaluating the two expressions at $z_1$, we obtain:
\begin{equation}
(\frac{df^2}{dz}(z_1)z_1 - f^2(z_1))(1-\frac{\kappa_1}{\kappa_2})=0,
\end{equation}
and using $2f + \frac{df}{dz}z=0$ we obtain $-5
f^2(z_1)(1-\frac{\kappa_1}{\kappa_2})=0$. For $\kappa_2=\kappa_1$, the
solution can be found~[2] and thus the assumption of complete
similarity in $\Pi_2$ is correct. However, in the case of
$\kappa_2\ne\kappa_1$, we have $f(z_1)=0$. This is a contradiction,
because the change in sign of $\der{h(x,t)}{t}$ should occur inside
the mound, where the height is positive. Hence, the assumption of
complete similarity for $\kappa_2\ne\kappa_1$ does not hold.

We next solve the problem numerically and study the asymptotic
behavior of the solution.

\section{Numerical solution of the partial differential equation and
further analysis for problem 1.}\label{secpde}

In order to simplify the numerical solution for equation
(\ref{eq1}) with free boundary conditions (\ref{eq2}), we use a change
of variables: $\xi=\frac{x}{x_r}$. We set $H(\xi,t)=h(x_r \xi,t)$, and equation (\ref{eq1}) is
transformed:
\begin{equation}\label{eq10.1}
  \partial_{t}H= \begin{cases} \frac{1}{x_r^2}(\kappa_1 \partial_{\xi
  \xi} H^2(\xi,t)+ \kappa_1 \xi \partial_\xi H(1,t) \partial_\xi
  H(\xi,t)) & \text{if $\partial_{t} H(\xi,t)>0$},\\
  \frac{1}{x_r^2}(\kappa_2 \partial_{\xi \xi} H^2(\xi,t)+ \kappa_1 \xi
  \partial_\xi H(1,t) \partial_\xi H(\xi,t)) & \text{if $\partial_{t}
  H(\xi,t)<0$}; \end{cases}
\end{equation}
with boundary conditions $H(0,t)=H(1,t)=0$. This effectively fixes the right
boundary at $\xi = 1$.

 The location of the free boundary $x_r(t)$ can be obtained in the
course of the numerical solution in the following way. We assume that the solution is
nearly stationary near the tip $x_r$ and $h(x,t) \approx
H(x-vt)$. Here, $v$ denotes the instantaneous speed of mound extension, which
changes slowly as a function of $t$. Then $\partial_{t}{h} \approx
-v\partial_{x}{h}$ near $x_r$.

Considering equation (\ref{eq1}) near the boundary $x_r$ of the mound,
where $h$ is small, we have: $$\partial_{t}{h} =\kappa_1\partial_{xx}
h^2(x,t) =2\kappa_1 ((\partial_{x}h(x,t))^2 + h \,
\partial_{xx} h(x,t)) \approx 2 \kappa_1 (\partial_{x}h(x,t))^2$$
so that
\begin{eqnarray}
  v(t)=-2\kappa_1 \partial_{x}h(x_r,t),~x_r(t)=\int_0^t v(t) \, dt
  +x_r(0).
\end{eqnarray}

\begin{figure}[htbp]
\centering \epsfig{figure=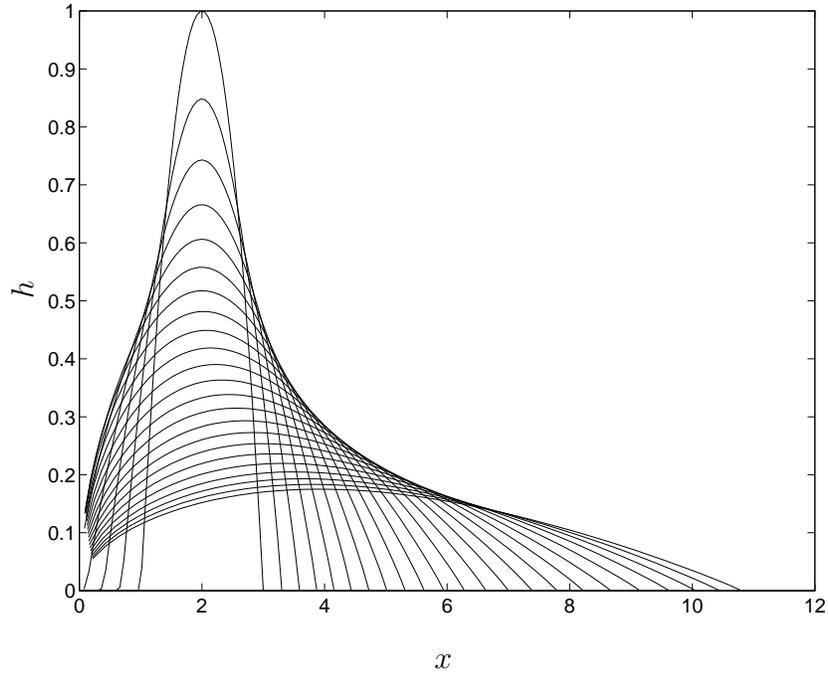,width=5in}
\put(-180,0){$x$}
\put(-340,140){{\begin{rotate} {$h$} \end{rotate}}}
\caption{Evolution in time of the partial differential equation
         solution with $\frac{\kappa_1}{\kappa_2}=.5$ for time
         interval t=[0,18]. } \label{fig:fig2}
\end{figure}

  We solve the new boundary value problem numerically by
using a forward-in-time, centered-in-space finite-difference
approximation, where $u_i^{n}$ is an approximation to the solution of
(\ref{eq10.1}) at the grid point $(x_i,t_n)$:
\begin{eqnarray}
 \kappa_i^n &=& \begin{cases} \kappa_1 & \text{if $ [(u_{i-1}^{n-1})^2
 -2(u_{i}^{n-1})^2 +(u_{i+1}^{n-1})^2] > 0$}, \\ \kappa_2 & \text{if $
 [(u_{i-1}^{n-1})^2 -2(u_{i}^{n-1})^2 +(u_{i+1}^{n-1})^2] < 0$}; \\
 \end{cases}\nonumber \\ u_i^{n+1}&=& u_i^{n}+\frac{\Delta t}{\Delta
 x^2} \{ \kappa_i^n [(u_{i-1}^n)^2 - 2(u_{i}^n)^2
 +(u_{i+1}^n)^2]\nonumber \\ &-& \kappa_1 \xi_i
 (u_i^n-u_{i-1}^n)(u_N^n-u_{N-1}^n)\}/{(x_r^n)^2}, \nonumber\\
 x_r^{n+1}&=& x_r^{n}-2 \kappa_1 \frac{\Delta t}{\Delta x} \>
 \frac{u_N-u_{N-1}}{x_r^{n}} . \nonumber
\end{eqnarray}
In the numerical computation we start with an initial distribution of
the source type, localized near $x=0$ (Fig.\ref{fig:fig2}). Before the left free boundary
reaches the point $x=0$, the solution is of the source type and we
can compare our numerical results to those in~[6]. After some
time $t$ the left free boundary reaches $x=0$, where it is thereafter
fixed (Fig. \ref{fig:fig2}).

\begin{figure}[htbp]
\centering \epsfig{figure=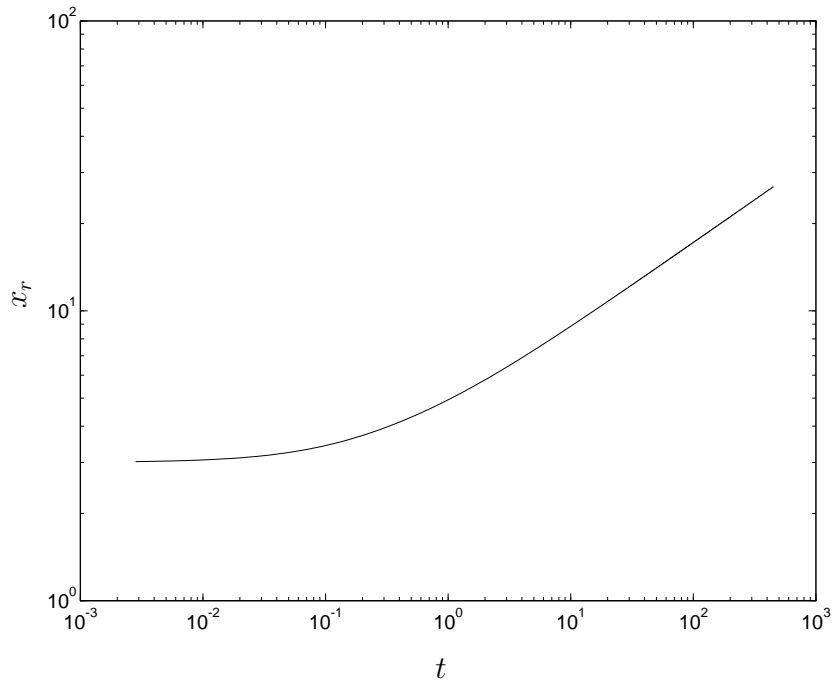,width=5in}
\put(-180,0){$t$}
\put(-340,140){{\begin{rotate} {$x_r$} \end{rotate}}}
\caption{Example of the time evolution for the position $x_r$ of the free boundary. Time interval $t=[0,450]$. For this plot $\frac{\kappa_1}{\kappa_2}=.5$.} \label{fig:fig3}
\end{figure}

\begin{figure}[htbp]
\centering \epsfig{figure=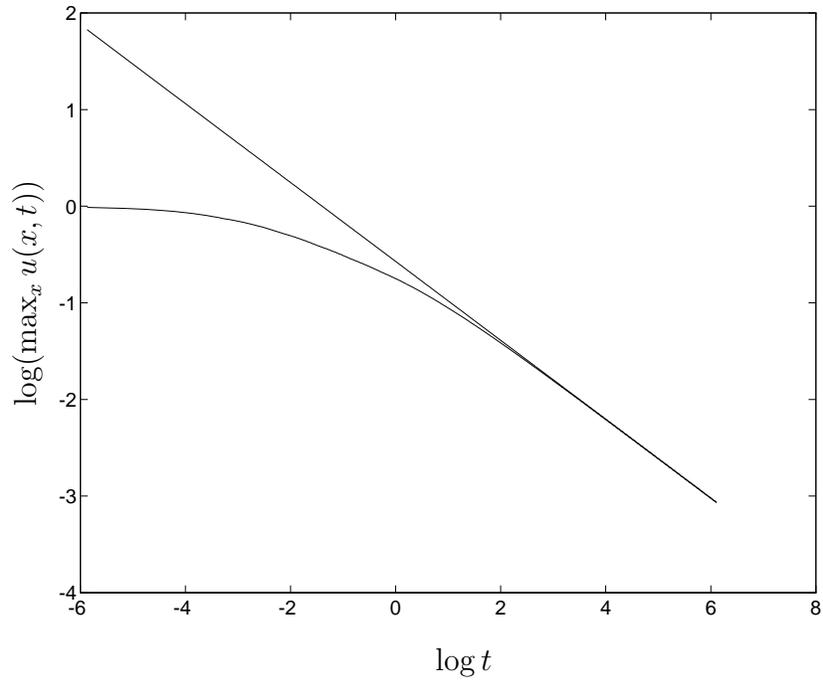,width=5in}
\put(-180,0){$\log t$}
\put(-340,100){{\begin{rotate} {$\log(\max_x u(x,t))$} \end{rotate}}}
\caption{Plot of $\log(\max_x u(x,t))$ vs. $\log(t)$ with $\frac{\kappa_1}{\kappa_2}=.5$. The straight line shows a linear fit to the straight part of the graph.} \label{fig:fig4}
\end{figure}

Now we consider the
scaled solution:
\begin{equation}\label{eq22}
\frac{H(\xi,t)}{\max_{\xi}\> H(\xi,t)} = f(\xi,\frac{\kappa_1}{\kappa_2}), ~\xi=\frac{x}{x_r}.
\end{equation} We can see in  Figs.
\ref{fig:fig6} and \ref{fig:fig7} that as time $t$ increases the
numerical solution approaches a self-similar regime, so that the graphs of
the scaled solution for different times ``collapse'' into a
single curve. Moreover, Figs. \ref{fig:fig3} and \ref{fig:fig4} show a
power-law dependence on time for both $\max_\xi H$ and $x_r$ in
the self-similar regime.

In part 3, we have shown that a self-similar solution of the first kind does
not exist for this problem. To explain what happened we return to the
dimensional analysis and look now for a generalized
self-similar solution. 

We have determined that the
variables in the problem are related as follows:
$\Pi=F(\Pi_1,\Pi_2,\Pi_3)$, where $$ \Pi_1=\frac{x}{(Q \kappa_1
t)^{1/4}},
\Pi_2=\frac{d}{(Q \kappa_1 t)^{1/4}},
\Pi_3=\frac{\kappa_1}{\kappa_2}, \text{and } \Pi=h(\frac{\kappa_1 t} {Q})^{1/2}.
$$
Our numerical investigation shows that for large $t$, as  $\Pi_2
\rightarrow 0$ :
\begin{equation}\label{eq8}
F(\Pi_1,\Pi_2,\Pi_3)=\Pi_2^{\gamma} f_2(\Pi_1\Pi_2^{-\eps},\Pi_3),
\end{equation}
where $\gamma$ and $\eps$ are constants. In fact, this is the next simplest
situation after complete self-similarity and it is referred to as
self-similarity of the second kind in $\Pi_2$ (see~[2]).

Indeed, from the analysis above:
\begin{eqnarray}\label{eq10}
\frac{\Pi_1}{\Pi_2^\eps} &=& x d^{-\eps} (t\kappa_1 Q)^\frac{\eps-1}{4},	\nonumber \\
\Pi_2^\gamma &=& \frac{d^\gamma} {(t \kappa_1 Q)^{\gamma/4}}, 			\nonumber \\
h(x,t) &=& A \; t^{-\alpha} f(\frac{x}{Bt^{\>\beta}},\frac{\kappa_1}{\kappa_2}),	\\
x_r(t)&=&Bt^\beta,								\nonumber
\end{eqnarray}
where $$A=(\kappa_1Q)^{-\gamma/4}d^\gamma,
B=d^\eps(\kappa_1Q)^\frac{\eps-1}{4},
\alpha=\frac{\gamma + 2}{4}, \text{and } \beta=\frac{1-\eps}{4}.$$

The parameters $\alpha$ and $\beta$ depend on the ratio
$\frac{\kappa_1}{\kappa_2}$. They cannot be determined on the basis of
dimensional analysis alone and have to be computed as a part of the
solution. We will see that there is, actually, only one unknown
parameter involved, since the differential equation provides an
additional relation between
$\alpha$ and $\beta$.

\section{Derivation and numerical solution of a nonlinear eigenvalue problem.}
The numerical solution of partial differential equation showed
that there is indeed an intermediate asymptotic solution of the form
(\ref{eq10}).  Now, we can obtain such a self-similar solution by
transforming the problem of
solving partial differential equation (\ref{eq1}) with boundary
conditions (\ref{eq2}) into a nonlinear eigenvalue problem.

 We substitute (\ref{eq10}) into (\ref{eq1}) and normalizing $f$
so that $\frac{A}{B^2}=1$ we get:
\begin{eqnarray} \label{eq11}
        \partial_t h &=& -t^{-(\alpha+1)} \; A(\alpha f+\beta f' \xi), \\
	 \partial_{xx} h^2 &=& \frac{A^2 t^{-2 \alpha} (f^2(\xi))''}{B^2 t^{2 \beta}},\\
	\alpha f(\xi) + \beta f'(\xi) \xi &=&
                -\kappa^* (f^2(\xi))'' t^{- \alpha -2\beta +1},
                \label{eq15}
\end{eqnarray}
where
\begin{equation}
\kappa^* = \begin{cases}
                \kappa_1 & \text{if $((1-2\beta)f+\beta f' \xi)<0$}, \\
                \kappa_2 & \text{if $((1-2\beta)f+\beta f' \xi)>0$}. \\
                \end{cases}\nonumber
\end{equation}
Since equation (\ref{eq15}) cannot depend on time explicitly, $ -
\alpha-2\beta +1=0$.  Finally, we get an ordinary differential
equation:
\begin{equation}\label{eq12}
(f^2)''= \begin{cases} -(\beta f'\xi +(1-2\beta) f) & \text{if
$((1-2\beta)f+\beta f' \xi) <0$}, \\ -\frac{\kappa_2}{\kappa_1}(\beta
f'\xi +(1-2\beta) f) & \text{if $((1-2\beta)f+\beta f' \xi)>0$}.
\end{cases}
\end{equation}
The boundaries $x_0=0$ and $x_r={B t^\beta}$, in the new space
variable, correspond to $\xi=0$ and $\xi=1$. The boundary condition at
$\xi=0$ becomes: 
\begin{equation}\label{eq14.1}
f(0)=0.
\end{equation}
For the right boundary, $\xi=1$, we have:
\begin{equation}\label{eq13}
f(1)=0 \>\text{and}\> (f^2)'(1)=0.
\end{equation}
From \ref{eq12} and \ref{eq13} it follows that $2\kappa_1
(f'(1))^2+\beta \xi f'(1)=0$ and the tip conditions become:
\begin{equation}\label{eq14}
f(1)=0 \>\text{and}\> f'(1)=-\frac{\beta}{2\kappa_1}.
\end{equation}

The second order ODE (\ref{eq12}) with three boundary conditions (\ref{eq14} and \ref{eq14.1}) constitutes a non-linear eigenvalue problem, which we now have to
solve numerically. For each value of $\frac{\kappa_1}{\kappa_2}$, we
find a value of $\beta$ such that the boundary conditions are satisfied.

\begin{figure}[htbp]
\centering \epsfig{figure=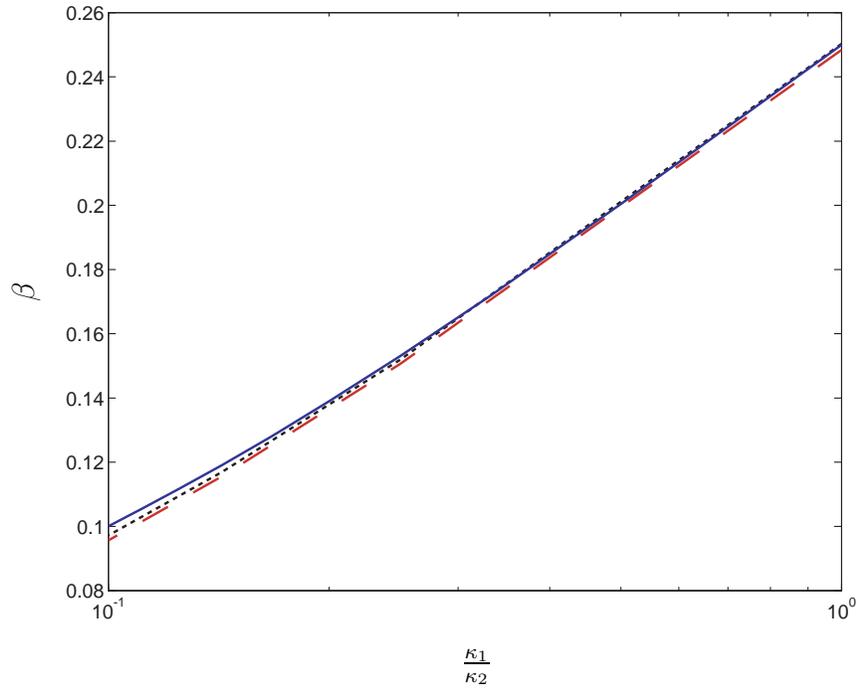,width=5in}
\put(-180,0){$\frac{\kappa_1}{\kappa_2}$}
\put(-350,140){{\begin{rotate} {$\beta$} \end{rotate}}}
\caption{Plot of the dependence of $\beta$ on $\frac{\kappa_1}{\kappa_2}$.
The solid line is obtained from graph \ref{fig:fig4} for different
values of $\frac{\kappa_1}{\kappa_2}$, by applying least-squares fit
to the straight part of the graph and then using equation
(\ref{eq20}). The dashed line is obtained from graphs \ref{fig:fig3},
by applying the procedure above and then (\ref{eq21}). The dotted line is
obtained from  the solution of the nonlinear eigenvalue problem given by
equation (\ref{eq12}) with initial conditions \ref{eq14} for different
values of $\frac{\kappa_1}{\kappa_2}$.}
\label{fig:fig5}
\end{figure}

\begin{figure}[htbp]
\centering \epsfig{figure=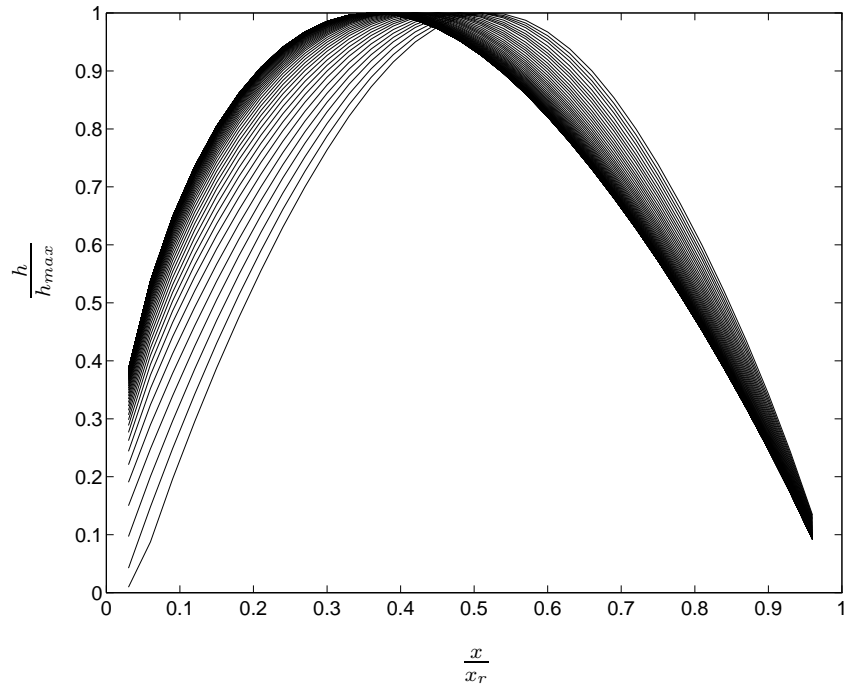,width=5in}
\put(-180,0){$\frac{x}{x_r}$}
\put(-350,140){{\begin{rotate} {$\frac{h}{h_{max}}$} \end{rotate}}}
\caption{Plot of the scaled solution (\ref{eq22}) of the partial differential equation for $t=[.2, 18]$. ($\frac{\kappa_1}{\kappa_2}=.5$).} \label{fig:fig6}
\end{figure}

\begin{figure}[htbp]
\centering \epsfig{figure=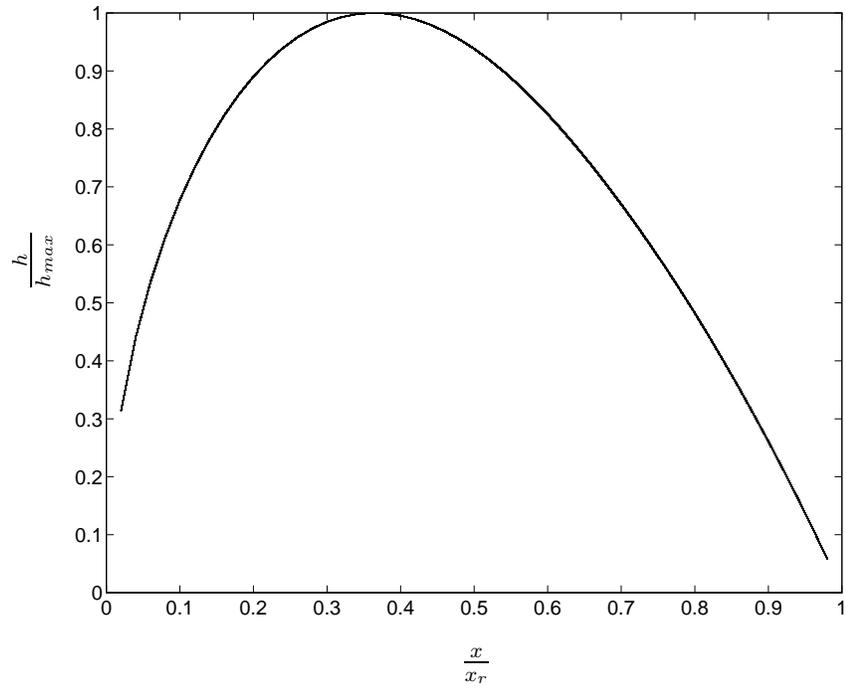,width=5in}
\put(-180,0){$\frac{x}{x_r}$}
\put(-350,140){{\begin{rotate} { $\frac{h}{h_{max}}$} \end{rotate}}}
\caption{Plot of the scaled solution (\ref{eq22}) of the partial differential equation for time interval $t=[18,450]$. For this plot $\frac{\kappa_1}{\kappa_2}=.5$.} \label{fig:fig7}
\end{figure}

We use a high order, Taylor-expansion-based method to start the
integration at $\xi=1$ followed by a 4th order Runge-Kutta method
and an iterative procedue to arrive at the value for $\beta$ such
that the third condition $f(0)=0$ is satisfied. For computational
convenience, we transform the differential equation by changing
variables: $g(\xi)=f^2(\xi)$, so that $g(\xi)$ does not have a
singularity at $\xi=0$. In this manner, we obtain the dependence of
$\beta$ on $\frac{\kappa_1}{\kappa_2}$ (Fig.
\ref{fig:fig5}).

\section{Comparison of the results for problem 1.}\label{sec_compar}

In logarithmic coordinates, we obtain from (\ref{eq10}):
\begin{eqnarray}
 \log(u(1/2,t))&=&-\alpha \log(t) + \log(A \>
 f(1/2,\frac{\kappa_1}{\kappa_2})) \label{eq20}\\ \log(x_r(t))&=&\beta
 \log(t)+\log(B) \label{eq21}
\end{eqnarray}
i.e., straight lines with slopes $-\alpha$ and $\beta$.

From plots in Figs. \ref{fig:fig3} and \ref{fig:fig4}, we can observe
that after some initial time $t$ both graphs approach straight
lines. We repeat the calculations for a range of values of
$\frac{\kappa_1}{\kappa_2}$.
   
Comparison of the results of the numerical solution of the nonlinear
eigenvalue problem with the results obtained from the numerical
solution to the partial differential equation (Fig. \ref{fig:fig5}),
shows that the two agree with high precision.

Also, the exact solution for the case $\kappa_1=\kappa_2$ gives the value
$\beta=.25$, which coincides with the results of the numerical
computations with good accuracy.

\section{Numerical solution of the partial differential equation for problem 2.}

Although problems 1 and 2 are similar, the numerical treatment of
problem 2 is more complicated. Time evolution of the left
boundary in problem 2 makes rescaling, which was used in the numerical
solution of problem 1, infeasible. Instead, we solve equation
(\ref{eq1}) on a grid, taking into account that the left and right
boundaries may not fall onto gridpoints. We determine new positions of
the boundaries from the numerical solution at each timestep.

 Equation (\ref{eq1}) is discretized using a forward-in-time, centered-in-space finite-difference scheme:
\begin{eqnarray}
 \kappa_i^n &= \begin{cases} \kappa_1 & \text{if } [(u_{i-1}^{n-1})^2
 -2(u_{i}^{n-1})^2 +(u_{i+1}^{n-1})^2] > 0, \\ \kappa_2 & \text{if }
 [(u_{i-1}^{n-1})^2 -2(u_{i}^{n-1})^2 +(u_{i+1}^{n-1})^2] < 0;\\
 \end{cases}\nonumber \\ 
u_i^{n+1}&= u_i^{n}+\frac{\Delta t}{\Delta x^2} 									\{ \kappa_i^n [(u_{i-1}^n)^2 - 2(u_{i}^n)^2										+(u_{i+1}^n)^2]\} ,		\nonumber \\
u_l^{n+1}&= u_l^{n}+\frac{2 \Delta t}{\Delta x + \Delta x_l} \kappa_2^n 						\{ \frac{(u_{l+1}^n)^2 - (u_{l}^n)^2}
					{\Delta x} - q^n \},				 \\ 
u_r^{n+1}&= u_r^{n}+\frac{2\Delta t}{\Delta x + \Delta x_r} \kappa_1^n 						\{ \frac{(u_{r-1}^n)^2 -(u_{r}^n)^2}{\Delta x} - 
					\frac{(u_{r}^n)^2}{\Delta x_r} \}.	\nonumber
\end{eqnarray}
Here $u_l$ and $u_r$ are the nonzero values of $u$ on the grid,
adjacent to the left and right boundaries respectively, $q$ is the
drainage flux, $\Delta x_l$ and $\Delta x_r$ are distances from the
left and right boundaries to the grid points. We treat the values
$u_l$ and $u_r$ separately in order to incorporate the boundary
conditions and improve precision.

The location of the left boundary is obtained from the values of
$u$:
\begin{equation}
x_l^{n+1}=x_l^n +\Delta x_l - \frac{(u_l^{n})^2}{q}.
\end{equation}
The right boundary location is obtained by extrapolation from
the values of $u^{n+1}$.

We check the numerical method for $\kappa_1=\kappa_2$ by comparing the
numerical solution with a known analytic solution. The exact self-similar solutions
for the problems with forced drainage are given in~[3]. We
choose a value of $\beta \ne .25$ and then solve an ordinary
differential equation (\ref{eq15}) with the initial condition
(\ref{eq14}). For $\beta < .25$ the solution $f(\xi)$ of the ordinary
differential equation intersects the $x$-axis at some point
$\xi=\lambda > 0$ and at $\xi=1$.

From the solution of the ordinary differential equation we obtain a
self-similar solution:
\begin{equation}
h(x,t)=B t^{-(1-2 \beta)} f(\frac{x}{A t^\beta})
\end{equation}
of the partial differential equation. The locations of the free
boundaries are given by $x_r(t)=A t^{\beta}$, $x_l(t)=\lambda A
t^{\beta}$, and the drainage flux is given by $q(t)=m B A t^{-2+3 \beta} (f^2)'(\lambda)$
(see~[3]). We use the self-similar solution at some time $t_0$ as
an initial value for the numerical solver, set drainage flux on the left
boundary to be $q(t)$, and compute the solutions until time
$t_1$. As in the analysis in section ~\ref{sec_compar}, the graphs
of $x_r(t)$, $x_l(t)$ should be straight lines in logarithmic
coordinates, and the graphs of the scaled solution for different times $t$
should collapse into one curve.  That's what we observe in Figs.
~\ref{fig:fig9} and ~\ref{fig:fig10}.

Now, we try to model the conditions of a flood followed by forced
drainage, as described in the introduction. We begin by computing the
solution to problem 1 until some time $t_0$, which corresponds to the
flood followed by natural drainage through the boundary of the
aquifer. After $t_0$, we set a constant drainage flux $q(t) = q_0$ at
the left boundary. In particular, we set $q_0$ to equal twice the
natural drainage flux at time $t_0$. As we see in Fig.
~\ref{fig:fig11}, the water mound, that has appeared after the flood, is
completely extinguished in finite time.

\section{Conclusion.}
\begin{enumerate}
\item The numerical simulations of two problems involving drainage and
capillary retention of the fluid a in porous medium were presented.  It
was shown that the problem with dipole type initial and boundary conditions has a
self-similar intermediate asymptotics in the case of a porous medium
with capillary retention.
\item A problem of control of the water mound
extension by forced drainage was considered. The possibility of
extinguishing the propagating water mound by creating a forced drainage
flux $q(t)$ at the left boundary was confirmed numerically. 
Using our results, it should be possible to derive a cost efficient drilling
regime and to localize the mound and contain the contamination inside
a prescribed region.
\end{enumerate}
It would be interesting to extend the numerical investigation above to
the case of a fissurized porous medium.

\section{Acknowledgements.}
The authors are grateful to Professor G.I. Barenblatt, without whose
direction and advice this work would not have been possible. The
authors use this occasion to thank Professor A. Chorin for many helpful
discussions of this work and for his constant attention and encouragement. 

This work was supported in part by the Computational Science Graduate
Fellowship Program of the Office of Scientific Computing in the
Department of Energy, NSF grant contract DMS-9732710, and the Office of
Advanced Scientific Computing Research,
Mathematical, Information, and Computational Sciences Division,
Applied Mathematical Sciences Subprogram, of
the U.S. Department of Energy, under Contract No. DE-AC03-76SF00098.

\section*{References}

\begin{itemize}
\item[1.]
G.I. Barenblatt, \emph{Scaling, self-similarity, and intermediate asymptotics},
  first ed., Cambridge University Press, New York, 1996.

\item[2.]
G.I. Barenblatt, V.M. Entov, and V.M. Ryzhik, \emph{Theory of fluid flows
  through natural rocks}, first ed., Kluwer Academic Publishers, Dordrecht,
  1990.

\item[3.]
G.I. Barenblatt and J.L. Vasquez, \emph{A new free boundary problem for
  unsteady flows in porous media}, Euro. Jnl of Applied Mathematics \textbf{9}
  (1998), 37--54.

\item[4.]
C.W. Fetter, \emph{Applied hydrogeology}, third ed., Macmillan College
  Publishing Company, New York, 1988.

\item[5.]
A.S. Kalashnikov, \emph{Some problems of qualitative theory of the non-linear
  second-order parabolic equations}, Russian Math. Surveys (1987), no.~42,
  169--222.

\item[6.]
I.N. Kochina, N.N. Mikhailov, and M.V. Filinov, \emph{Groundwater mound
  damping.}, Int. J. Engng Sci \textbf{21} (1983), no.~4, 413--421.

\item[7.]
B.A. Wagner, \emph{Perturbation techniques and similarity analysis for the
  evolution of interfaces in diffusion and surface tension driven problems},
  Zentrum Mathematik, TU Munchen, 1999.
\end{itemize}
\newpage
\begin{figure}[htbp]
\centering \epsfig{figure=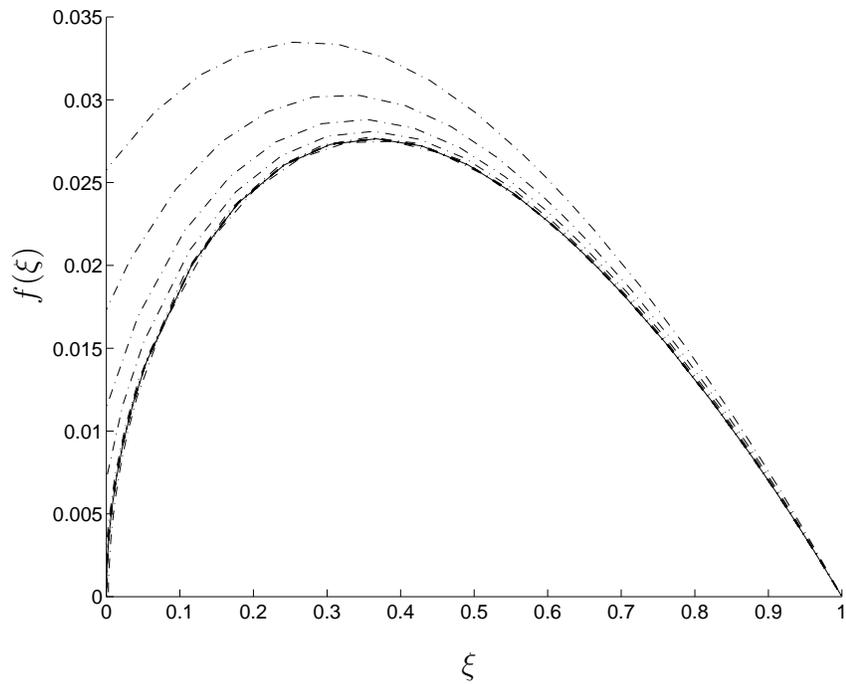,width=5in}
\put(-180,0){$\xi$}
\put(-350,140){{\begin{rotate} {$f(\xi)$} \end{rotate}}}
\caption{Example of convergence of the ``shooting'' method for solving the nonlinear eigenvalue problem given by (\ref{eq12}), (\ref{eq14}).} \label{fig:fig8}
\end{figure}
\begin{figure}[htbp]
\centering \epsfig{figure=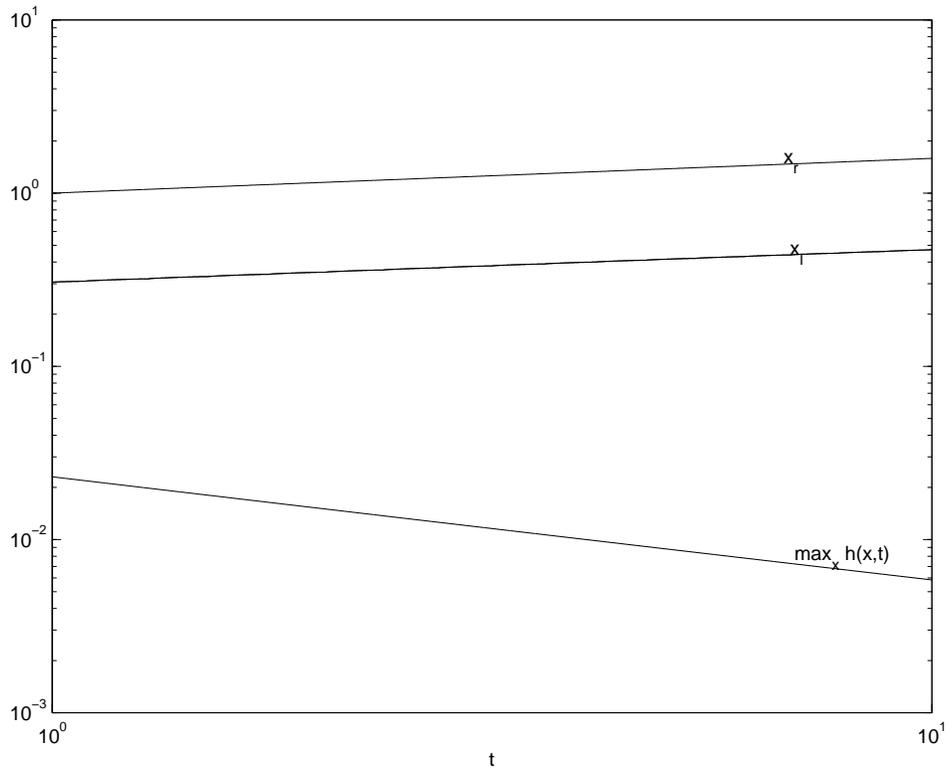,width=5in}
\caption{Positions of the boundaries and $\max_x h(x,t)$ vs $t$ in logarithmic coordinates for the self-similar regime in the forced drainage problem. $\beta=.2$} \label{fig:fig9}
\end{figure}
\begin{figure}[htbp]
\centering \epsfig{figure=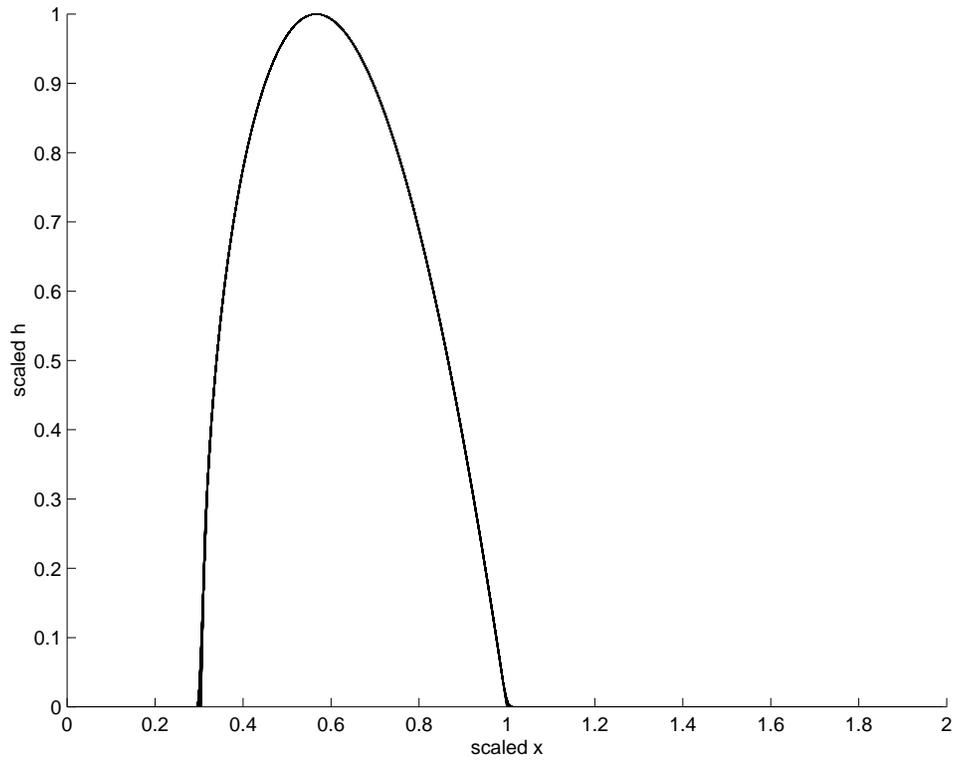,width=5in}
\caption{``Collapse'' of the graphs of the scaled solution in self-similar coordinates for problem 2. Time $t=[1,10]$, $\beta=.2$ } \label{fig:fig10}
\end{figure}
\begin{figure}[htbp]
\centering \epsfig{figure=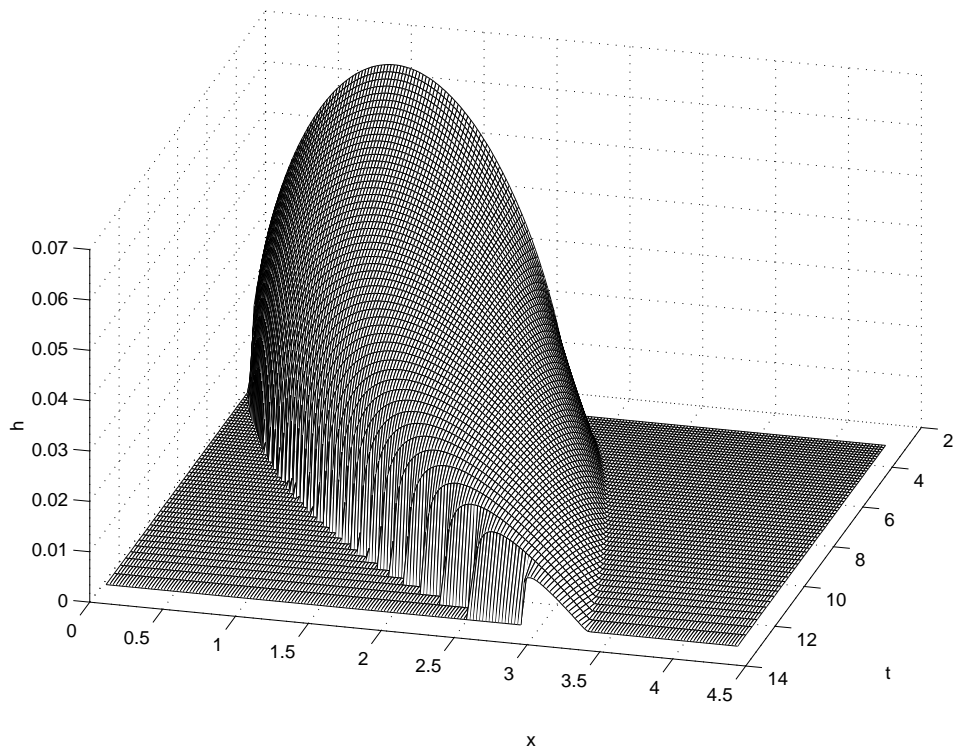,width=5in}
\caption{Evolution in time of the solution to the PDE with forced drainage.}
\end{figure}
\begin{figure}[htbp]
\centering \epsfig{figure=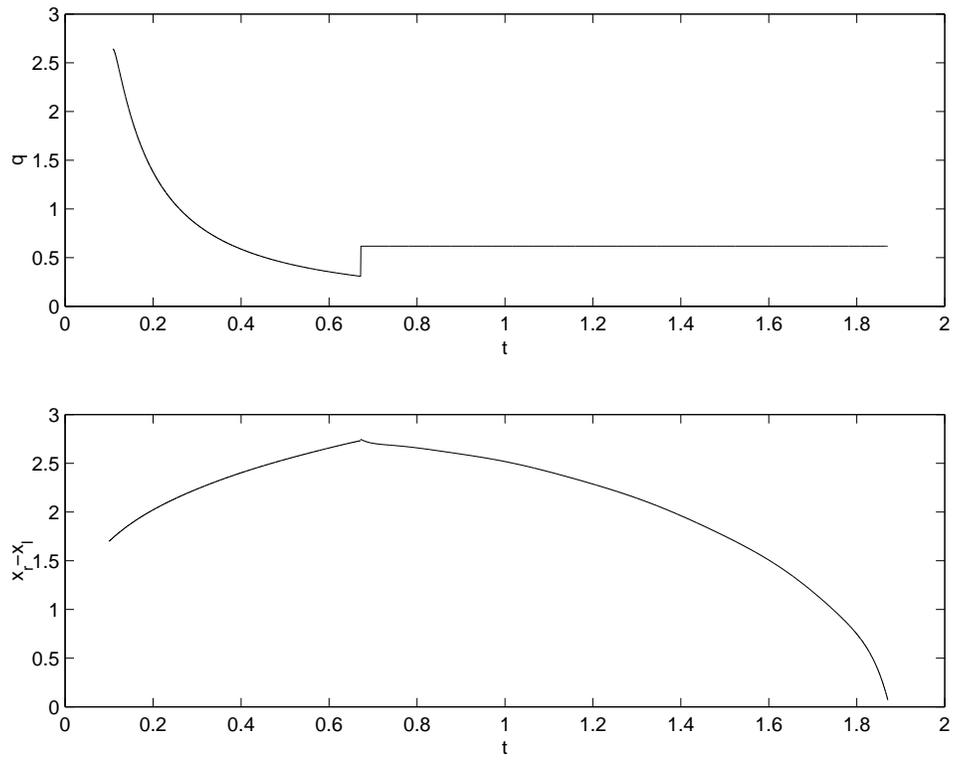,width=5in}
\caption{Extinguishing of the water mound with forced drainage.} \label{fig:fig11}
\end{figure}
\end{document}